\documentclass[11pt,a4paper]{article}

\usepackage{graphicx}        
\usepackage[bottom]{footmisc}
\usepackage{url}
\usepackage{xcolor}
\usepackage{amsmath,amsthm,amssymb}
\usepackage{mathtools}
\usepackage{authblk}
\usepackage[labelfont=bf]{caption}
\usepackage{hyperref}
\usepackage{geometry}

\newtheorem{theorem}{Theorem}
\newtheorem{lemma}{Lemma}

\begin{document}
\title{Space-time hexahedral finite element methods for parabolic evolution problems\thanks{This work was supported by the Austrian Science Fund (FWF) under grant W1214, project DK4.}}
\author[1]{Ulrich Langer}
\affil[1]{%
    Institute for Computational Mathematics\authorcr
    Johannes Kepler University Linz\authorcr
    Altenberger Str. 69, 4040 Linz, Austria}
\author[2]{Andreas Schafelner}
\affil[2]{%
    Doctoral Program ``Computational Mathematics''\authorcr
    Johannes Kepler University Linz\authorcr
    Altenberger Str. 69, A-4040 Linz, Austria}
\date{}
%
%
\maketitle

\abstract
{We present locally stabilized, conforming space-time finite element methods
    for parabolic evolution equations
    on hexahedral decompositions of the space-time cylinder.
    Tensor-product decompositions allow for anisotropic a priori error estimates, that are explicit in spatial and temporal meshsizes.
    Moreover, tensor-product finite elements are suitable for anisotropic adaptive mesh refinement strategies
    provided that an appropriate a posteriori discretization error estimator is available.
    We present such anisotropic adaptive strategies together with numerical experiments.\\[1em]
}

%
%
\section{Introduction}
\label{LS:sec:Introduction}
We
consider the parabolic initial-boundary value problem (IBVP),
find  $ u $ such that
\begin{equation}
    \label{LS:eqn:IBVP}
    \partial_{t}u - \mathrm{div}_x(\nu\, \nabla_{x} u ) = f\ \text{in}\ Q,\;
    u = u_D := 0\ \text{on}\ \Sigma,\;
    u = u_0 := 0\ \text{on}\ \Sigma_0,\;
\end{equation}
as a model problem typically arising in heat conduction and diffusion,
where $Q=\Omega\times(0,T)$, $\Sigma = \partial\Omega\times(0,T)$, and
$\Sigma_0 = \Omega\times\{0\}$.
The spatial domain $\Omega \subset \mathbb{R}^d$, $d=1,2,3$, is assumed to be bounded and Lipschitz,
$T > 0$ is the terminal time, $ f \in L_2(Q) $ denotes a given source,
and
$\nu \in L_\infty(Q)$ is a given uniformly bounded and positive coefficient
that may discontinuously depend on the spatial variable $x =(x_1,\ldots,x_d)$ and the time variable $t$,
but $\nu(x,t)$ should be of bounded variation in $t$ for almost all $x \in \Omega$.
Then there is a unique weak solution
$u \in V_0 := \{v \in  L_2(0,T;H^1_0(\Omega)): \partial_t u \in  L_2(0,T;H^{-1}(\Omega)), v=0 \mbox{ on } \Sigma_0  \} $
of the IBVP \eqref{LS:eqn:IBVP}; see, e.g.,
\cite{LS:Evans:2010a,LS:Steinbach:2015a}.
Moreover,  $\partial_t u$ and $Lu := - \mathrm{div}_x(\nu\, \nabla_{x} u )$ belong to $L_2(Q)$;
see \cite{LS:Dier:2015a}.
The latter property is called maximal parabolic regularity.
In this case, the parabolic partial differential equation
$\partial_{t}u - \mathrm{div}_x(\nu\, \nabla_{x} u ) = f$
holds in $L_2(Q)$.
This remains even valid for inhomogeneous initial conditions $u_0 \in H^1_0(\Omega)$.

Time-stepping methods in combination with some
spatial
discretization method
like the finite element method (FEM) are still the standard approach
to the numerical solution of IBVP like \eqref{LS:eqn:IBVP};
see, e.g., \cite{LS:Thomee2006a}. This time-stepping approach as well as the more recent
discontinuous Galerkin, or
discontinuous Petrov-Galerkin
methods based on time slices or slabs are in principle sequential.
The sequential nature of these methods hampers the full space-time adaptivity
and parallelization; but see the overview paper \cite{LS:Gander2015a}
for parallel-in-time methods. Space-time finite element methods on fully unstructured
decomposition of the space-time cylinder $Q$ avoid these bottlenecks;
see \cite{LS:SteinbachYang:2019a} for an overview of such kind of space-time methods.

In this paper, we follow
our preceding papers  \cite{LS:LangerNeumuellerSchafelner:2019a,LS:LangerSchafelner:2020a,LS:LangerSchafelner:2020b},
and construct locally stabilized, conforming space-time finite element schemes
for solving the IBVP \eqref{LS:eqn:IBVP},
but on hexahedral meshes that are more suited for anisotropic refinement
than simplicial meshes used in
\cite{LS:LangerNeumuellerSchafelner:2019a,LS:LangerSchafelner:2020a,LS:LangerSchafelner:2020b}.
We mention that SUPG/SD and Galerkin/least-squares
stabilizations of time-slice finite element schemes for solving  transient
problems were already used in early papers; see, e.g.,
\cite{LS:JohnsonSaranen:1986a} and
\cite{LS:HughesFrancaHulbert:1989a}.
Section~\ref{LS:sec:SpaceTimeFEM} recalls the construction
of locally stabilized space-time finite element schemes, the properties
of the corresponding discrete bilinear form, and the a priori
discretization error estimates from
\cite{LS:LangerNeumuellerSchafelner:2019a,LS:LangerSchafelner:2020a,LS:LangerSchafelner:2020b}.
In Section~\ref{LS:sec:ErrorEstimates}, we derive new anisotropic a priori discretization estimates
for hexahedral tensor-product meshes, and we provide anisotropic adaptive mesh refinement strategies
that are based on a posteriori error estimates, anisotropy indicators, and anisotropic
adaptive mesh refinement using hanging nodes.
In Section~\ref{LS:sec:NumericalResults}, we present and discuss numerical results
for
an example where a singularity occurs in the spatial gradient of the solution.
The large-scale system of space-time finite element equations is always solved by means of
the Flexible Generalized Minimal Residual (FGMRES) method preconditioned by space-time algebraic multigrid.
%
%
%
\section{Space-time finite element methods}
\label{LS:sec:SpaceTimeFEM}
In this section, we will briefly describe the
space-time finite element method based on localized time-upwind stabilizations;
for details of the construction and analysis, we refer to our previous
work \cite{LS:LangerNeumuellerSchafelner:2019a,LS:LangerSchafelner:2020a,LS:LangerSchafelner:2020b}.
Let $\mathcal{T}_h$ be a shape regular decomposition of the space-time cylinder $Q$, i.e., $\overline{Q} = \bigcup_{K\in\mathcal{T}_h} \overline{K}$, and  $K\cap K'=\emptyset$
for all $K$ and $K'$ from $ \mathcal{T}_h $ with $K \neq K'$;
see, e.g.,
\cite{LS:Ciarlet:1978a}  for more details.
Furthermore, we assume that $\nu$ is piecewise smooth, and possible
discontinuities are aligned with the triangulation as usual.
On the basis of the triangulation $\mathcal{T}_h$, we define the space-time finite element space
\begin{equation*}
    V_{0h} = \{ v\in C(\overline{Q})  : v(x_K(\cdot)) \in\mathbb{P}_p(\hat{K}),\,
    \forall  K \in \mathcal{T}_h,\, v=0\;\mbox{on}\; {\overline \Sigma} { \cap {\overline \Sigma}_0 }
    \},
\end{equation*}
where $x_K(\cdot)$ denotes the map from the reference element $\hat{K}$ to the finite element
$K \in \mathcal{T}_h$,
and $\mathbb{P}_p(\hat{K})$ is the space of polynomials of the degree $p$ on the reference element $\hat{K}$.
Since we are in the maximal parabolic regularity setting,
the parabolic Partial Differential Equation  (PDE) is valid in $L_2(Q)$.
Multiplying the PDE \eqref{LS:eqn:IBVP}, restricted to  $K \in\mathcal{T}_h$,
by a locally scaled upwind test function
$
    v_{h,K}(x,t) := v_h(x,t) + \theta_K h_K \partial_{t} v_h(x,t),\; v_h \in  V_{0h},
$
integrating over $K$, summing over all elements, applying integration by parts, and incorporating the
Dirichlet boundary conditions,
we obtain the variational {consistency} identity
\begin{equation}\label{LS:eqn:VariationalIdentity}
    a_h(u,v_h) = \ell_h(v_h),\quad \forall v_h\in V_{0h},
\end{equation}
with the mesh-dependent bilinear form
\begin{align}
    a_h(u,v_h) = & \sum_{K\in\mathcal{T}_h} \int_{K} \big[ \partial_t u v_h + \theta_K h_K \partial_t u \partial_t v_h \label{LS:eqn:bilinearform} \\ \nonumber
                 & \qquad\qquad + \nu \nabla_{x}u\cdot\nabla_{x}v_h - \theta_K h_K \mathrm{div}_x(\nu \nabla_{x} u)\partial_t v_h \big]
    \mathrm{d}K, \nonumber
\end{align}
and the mesh-dependent linear form
\begin{align*}
    \ell_h(v_h) = & \sum_{K\in\mathcal{T}_h} \int_{K} \big[f v_h + \theta_K h_K f \partial_t v_h\big]
    \mathrm{d}K.
\end{align*}
Now we apply the Galerkin principle, i.e., we look for a finite element approximation $u_h \in V_{0h}$ to $u$ such that
\begin{equation}\label{LS:eqn:SpaceTimeFiniteElementScheme}
    a_h(u_h,v_h) = l_h(v_h),\quad \forall v_h\in V_{0h}.
\end{equation}
Using Galerkin orthogonality (subtracting \eqref{LS:eqn:SpaceTimeFiniteElementScheme}  from   \eqref{LS:eqn:VariationalIdentity}),
and coercivity and extended boundedness of the bilinear form \eqref{LS:eqn:bilinearform}, we can show the following C\'{e}a-like best approximation estimate;
see \cite{LS:LangerNeumuellerSchafelner:2019a,LS:LangerSchafelner:2020a,LS:LangerSchafelner:2020b} for the proofs.
\begin{theorem}\label{LS:thm:Cea}
    Let
    $u \in H^{L,1}_{0}(Q) := \{v \in V_0\cap H^1(Q): Lv := -\mathrm{div}_x (\nu \nabla_{x} v) \in L_2(Q)\}$
    and $u_h \in V_{0h}$ be the solutions of the parabolic IBVP \eqref{LS:eqn:IBVP}     and the space-time finite element scheme \eqref{LS:eqn:SpaceTimeFiniteElementScheme}, respectively.
    Then the discretization error estimate
    \begin{equation}
        \label{LS:eqn:Cea}
        \|u-u_h\|_h \leq \inf_{v_h\in V_{0h}} \left(\|u-v_h\|_h + \frac{\mu_b}{\mu_c}\|u-v_h\|_{h,*}\right)
    \end{equation}
    is valid provided that $\theta_K = O(h_K)$ is sufficiently small, where
    \begin{align*}
        \| v \|_{h}^2   & = \frac{1}{2}\| v \|_{L_2(\Sigma_T)}^2 + \sum_{K\in\mathcal{T}_h} \Bigl[ \theta_{K} h_K\| \partial_{t} v \|_{L_2(K)}^2 + \|\nu^{1/2} \nabla_x v \|_{L_2(K)}^2 \Bigr], \\
        \| v \|_{h,*}^2 & =
        \|v\|_{h}^2 + \sum_{K\in\mathcal{T}_h} \Bigl[  (\theta_{K} h_K)^{-1}\| v \|_{L_2(K)}^2 +
        \theta_{K} h_K \|\mathrm{div}_x (\nu \nabla_{x} v) \|_{L_2(K)}^2 \Bigr].
    \end{align*}
\end{theorem}
The best-approximation error estimate \eqref{LS:eqn:Cea} now leads to convergence
rate estimates under additional regularity assumptions.
If the solution $u$ of \eqref{LS:eqn:IBVP} belongs to
$H^{L,1}_{0}(Q) \cap H^l(Q)$, $l > 1$,
then $\|u-u_h\|_h \le c(u) h^{s-1}$, where $s = \min \{l,p+1\}$, $h = \min_{K\in\mathcal{T}_h} h_K$,
$c(u)$ depends on the regularity of $u$. We refer the reader to
\cite[Theorem 13.3]{LS:LangerNeumuellerSchafelner:2019a} and
\cite[Theorem 3]{LS:LangerSchafelner:2020a}
for the proof of more detailed estimates in terms of the local mesh-sizes $h_k$
and the local regularity of the solution $u$.
%
%
\section{Anisotropic a priori and a posteriori error estimates}
\label{LS:sec:ErrorEstimates}
The convergence rate estimates presented at the end of the previous section
consider only isotropic finite elements, but in many application the
solution $u$ evolves differently with respect to time and space directions.
So, we should permit anisotropic finite elements with different mesh sizes
in different directions.
This raises the question whether we can obtain (localized) a priori estimates that are explicit in spatial and temporal mesh sizes as well as in spatial and temporal regularity assumptions imposed on the solution $u$.
We refer to \cite{LS:Apel:1999a} for a comprehensive summary of anisotropic finite elements.
For the remainder of this section, we will now assume that $K$ is a \emph{brick element}
(hexahedral element for  the case $d=2$),
i.e., the edges of $K$ are parallel to the coordinate axes. Moreover, we assume that $ u \in H^{l}(Q)$ with $ l > (d+1)/2 $ an integer such that we can use the Lagrange interpolation operator $ I_h $.
Let $ h_{K,i} = \max\{ |x_i-x_i'| : x,x'\in K\}$, then $ h_{K,x} = \max_{i=1,\dots,d} h_{K,i} $ and $ h_{K,t} = h_{K,d+1}$. Furthermore, let $ e_h = u - I_h u $ and $ s = \min\{l, p+1\} $, where $p$ is the polynomial degree of the finite element
shape functions in every coordinate direction.
Using the anisotropic interpolation error estimates from \cite{LS:Apel:1999a},
we get
\begin{gather*}
    \|e_h\|_{L_2(K)}^2 \le c \left( \sum_{j=1}^{d} h_{K,j}^{2s} \|\partial_{x_j}^{s} u\|_{L_2(K)}^2 + h_{K,t}^{2s} \|\partial_{t}^{s} u\|_{L_2(K)}^2 \right),\\
    \|\partial_{t}(e_h)\|_{L_2(K)}^2 \le
    c \left( \sum_{j=1}^{d} h_{K,j}^{2(s-1)} \|\partial_{t}\partial_{x_j}^{(s-1)} u\|_{L_2(K)}^2 + h_{K,t}^{2(s-1)} \|\partial_{t}^{s} u\|_{L_2(K)}^2 \right),\\
    \|\partial_{x_i}(e_h)\|_{L_2(K)}^2 \le
    c \left( \sum_{j=1}^{d} h_{K,j}^{2(s-1)} \|\partial_{x_i}\partial_{x_j}^{(s-1)} u\|_{L_2(K)}^2 + h_{K,t}^{2(s-1)} \|\partial_{x_i}\partial_{t}^{(s-1)} u\|_{L_2(K)}^2 \right),\\
    \|\partial_{x_i}^2(e_h)\|_{L_2(K)}^2 \le
    c \left( \sum_{j=1}^{d} h_{K,j}^{2(s-2)} \|\partial_{x_i}^2\partial_{x_j}^{(s-2)} u\|_{L_2(K)}^2 + h_{K,t}^{2(s-2)} \|\partial_{x_i}^2\partial_{t}^{(s-2)} u\|_{L_2(K)}^2 \right),
\end{gather*}
for $ i = 1,\dots,d$,
where $c$ denotes generic positive constants.
These estimates of the interpolation error and its derivatives immediately lead
to the corresponding interpolation error estimates with respect to the norms
$\|\cdot\|_{h}$ and $\|\cdot\|_{h,*}$.
\begin{lemma}\label{LS:lem:AnisotropicInterpolationErrorEstimates}
    Let $ u \in H^{L,1}_0(Q) \cap H^l(Q) $, $ l \in \mathbb{N} $ with $ l > (d+1)/2 $, and let $ \mathcal{T}_h $ be a decomposition of $Q$ into {brick elements}. Then the anisotropic interpolation error estimates
    \begin{align}
        \|u-I_h u\|_{h}   & \le \left(\sum_{K\in\mathcal{T}_h} h_{K,x}^{2(s-1)} \mathfrak{c}_{1}(u,K) + h_{K,t}^{2(s-1)} \mathfrak{c}_{2}(u,K)\right)^{1/2}, \label{LS:eq:AnisotropicInterpolationError:h}     \\
        \|u-I_h u\|_{h,*} & \le \left(\sum_{K\in\mathcal{T}_h} h_{K,x}^{2(s-1)} \mathfrak{c}_{1,*}(u,K) + h_{K,t}^{2(s-1)} \mathfrak{c}_{2,*}(u,K)\right)^{1/2} \label{LS:eq:AnisotropicInterpolationError:h*}
    \end{align}
    hold, where $ s = \min\{l, p+1\} $,
    and $\mathfrak{c}_{1}(u,K)$, $\mathfrak{c}_{2}(u,K)$, $\mathfrak{c}_{1,*}(u,K)$ and $\mathfrak{c}_{2,*}(u,K)$
    can easily be computed from the interpolation error estimates given above.
\end{lemma}
Combining the interpolation error estimates
\eqref{LS:eq:AnisotropicInterpolationError:h} and \eqref{LS:eq:AnisotropicInterpolationError:h*}
with the best approximation estimate \eqref{LS:eqn:Cea},
we can immediately derive an a priori discretization error estimate.
\begin{theorem}
    Let the assumptions of Theorem~\ref{LS:thm:Cea} (best approximation estimate) and of Lemma~\ref{LS:lem:AnisotropicInterpolationErrorEstimates} (anisotropic interpolation error estimates) be fulfilled. Then
    the anisotropic  a priori discretization error estimate
    \begin{equation*}
        \|u-u_h\|_{h} \le \left(\sum_{K\in\mathcal{T}_h} h_{K,x}^{2(s-1)} \mathfrak{C}_{1}(u,K) + h_{K,t}^{2(s-1)} \mathfrak{C}_{2}(u,K)\right)^{1/2},
    \end{equation*}
    is valid, where $ s = \min\{l, p+1\} $,
    and $\mathfrak{C}_{1}(u,K)$ and $\mathfrak{C}_{2}(u,K)$ can be computed from
    \eqref{LS:eq:AnisotropicInterpolationError:h} and \eqref{LS:eq:AnisotropicInterpolationError:h*}.
\end{theorem}
In the computational practice, we would like to replace the uniform mesh refinement
by adaptive space-time mesh refinement that takes care of possible anisotropic features
of the solution in space and time. Here brick finite elements with hanging nodes,
as implemented in MFEM (see next section), are especially suited.
To drive anisotropic adaptive mesh refinement, we need a localizable a posteriori
error estimator providing local error indicators, and an anisotropy indicator
defining the refinement directions in each brick element $K \in \mathcal{T}_h$.

We use the functional a posteriori error estimators introduced by Repin; see his
monograph \cite{LS:Repin:2008a}. 
Repin proposed two error majorants
$\mathfrak{M}_1$ and $\mathfrak{M}_2$
from which the local error indicators
    {%
        \[
            \eta_{1,K}^2(u_h) = \frac{1}{\delta}\int_{K}\!\big[ (1+\beta) |\mathbf{y} -\nu\nabla_{x}u_h|^2
                + \frac{1+\beta}{\beta} c_{F\Omega}^2 |f - \partial_tu_h + \mathrm{div}_x \mathbf{y}|^2\big]\mathrm{d}K
            \;\mbox{and}
        \]
        %
        \begin{align*}
            \eta_{2,K}^2(u_h) = & \frac{1}{\delta}\int_{K}\! (1+\beta)\bigl[ |\mathbf{y}-\nu\nabla_{x}u_h + \nabla_x \vartheta|^2 + \frac{c_{F\Omega}^2(1+\beta)}{\beta}  |f - \partial_tu_h - \partial_t \vartheta + \mathrm{div}_x \mathbf{y}|^2\bigr]\mathrm{d}K \\
                                & + \gamma \|\vartheta(\cdot,T)\|_{\Omega}^2
            + 2 \int_K\!\big[ \nabla_{x}u_h\cdot\nabla_x\vartheta + (\partial_tu_h - f)\vartheta \big] \mathrm{d}K
        \end{align*}
    }
can be derived for each element $K \in \mathcal{T}_h$,
where $ \mathbf{y} \in H(\mathrm{div}_x,Q)$ is an arbitrary approximation
to the flux, $\vartheta \in H^1(Q)$ is also an arbitrary function,
$\delta \in (0,2]$, $ \beta > \mu$, $\mu \in (0,1) $,  and $ \gamma > 1 $.
The positive constant $c_{F\Omega}$ denotes the constant in the inequality
$\|v\|_{L_2(Q)} \le c_{F\Omega} \|\sqrt{\nu} \nabla_x v\|_{L_2(Q)}$ for all $v \in L_2(0,T; H_0^1(\Omega))$,
which is nothing but the Friedrichs constant for the spatial domain $\Omega$ in the case $\nu = 1$.
Both majorants provide a guaranteed upper bound for the errors
\begin{gather*}
    |\!|\!|u-u_h|\!|\!|_{(1,2-\delta)}^2 \le \sum_{K\in\mathcal{T}_h}  \eta_{1,K}^2(u_h) \quad\text{and}\quad
    |\!|\!|u-u_h|\!|\!|_{(1-\frac{1}{\gamma},2-\delta)}^2 \le \sum_{K\in\mathcal{T}_h} \eta_{2,K}^2(u_h),
\end{gather*}
where
$|\!|\!|v|\!|\!|_{(\epsilon,\kappa)}^2 \coloneq \kappa \| \sqrt{\nu} \nabla_{x} v \|_{Q}^2 + \epsilon \|v\|_{\Sigma_T}^2 $.
%
%
Once we have computed the local error indicators $\eta_K(u_h)$ for all elements $ K \in \mathcal{T}_h $, we use \emph{D\"{o}rfler marking} to determine a set $ \mathcal{M} \subseteq \mathcal{T}_h $ of elements that will be marked for refinement. The set $ \mathcal{M} $ is of (almost) minimal cardinality such that
\[ \sigma\ \sum_{K\in\mathcal{T}_h} \eta_K(u_h)^2 \leq \sum_{K\in\mathcal{M}} \eta_K(u_h)^2, \]
where $ \sigma \in (0,1) $ is a bulk parameter. Let $ \mathbf{E}_K \in \mathbb{R}^{d+1}$ with entries $ E_i^{(K)}$, $i=1,\dots,d+1$, and $ \chi \in (0,1)$. In order to determine how to subdivide a marked element, we use the following heuristics: for each $ K \in \mathcal{M}$,
subdivide $K$ in direction $ x_i$ iff $E_i^{(K)}> \chi |\mathbf{E}_K|$.
In particular, we choose
\begin{equation*}
    \left(E_i^{(K)}\right)^2 := \begin{cases}
        \int_K\!({y}_{i} - \nu\,\partial_{x_i}u_h)^2\;\mathrm{d}K, & i\le d, \\
        \int_K\!(\mathbf{d_t}_h - \partial_{t}u_h)^2\;\mathrm{d}K, & i=d+1,
    \end{cases}
\end{equation*}
as our local anisotropy vector $ \mathbf{E}_K$, where $ \mathbf{y}_h = (y_i)_{i=1}^d$, $ \mathbf{d_t}_h = R_h(\partial_t u_h)$, and $ R_h $ is some nodal averaging operator like in a Zienkiewicz-Zhu approach.
%
%
\section{Numerical Results}
\label{LS:sec:NumericalResults}
Now let
$ \{p^{(j)}:j=1,\ldots,N_h\}$
be the finite element nodal basis of $ V_{0h}$,
i.e., $V_{0h} = \mbox{span}\{p^{(1)},\ldots,p^{(N_h)} \}$,
where $N_h$ is the number of all space-time unknowns (dofs).
Then we can express the approximate solution $ u_h $ in terms of this basis,
i.e., $ u_{h}(x,t) = \sum_{j=1}^{N_h} u_j\,p^{(j)}(x,t) $.
Inserting this
representation
into \eqref{LS:eqn:SpaceTimeFiniteElementScheme}, and testing with $p^{(i)}$,
we get the linear system
$
    {K}_{h} \underline{u}_{h}= \underline{f}_{h}
$
for determining the unknown coefficient vector $\underline{u}_{h}= (u_j)_{j=1,\ldots,N_h} \in \mathbb{R}^{N_h} $,
where $ {K}_{h} = (a_h(p^{(j)},p^{(i)}))_{i,j=1,\ldots,N_h} $ and
$ \underline{f}_{h} = (\ell_h(p^{(i)}))_{i=1,\ldots,N_h}$.
The system matrix $ {K}_{h}$ is non-symmetric, but positive definite due to coercivity
of the bilinear form $a_h(\cdot,\cdot)$. Thus, in order to obtain a numerical solution to the
to the IBVP \eqref{LS:eqn:IBVP}, we just need to solve one linear system of algebraic equations.
This is always solved by means of
the FGMRES method preconditioned by space-time algebraic multigrid (AMG). We use the finite element library MFEM \cite{LS:mfem-library} to implement our space-time finite element solver.
The AMG preconditioner is realized via \emph{BoomerAMG}, provided by the linear solver library \emph{hypre}\footnote{\url{https://github.com/hypre-space/hypre}}.
We start the linear solver with initial guess $ \mathbf{0} $, and stop once the initial residual has been reduced by a factor of $ 10^{-8}$.
In order to accelerate the solver in case of adaptive refinements, we also employ \emph{Nested Iterations}.
Here, we interpolate the finite element approximation from the previous
mesh to the current mesh, and use that as an initial guess for FGMRES. Moreover, we stop the linear solver earlier, e.g. once the residual is reduced by a factor of $ 10^{-2}$.
Furthermore, we will use the notation $ \mathcal{O}(h^{\alpha}) = \mathcal{O}(N_h^{-\alpha/(d+1)}) $ to indicate the corresponding convergence rates.

Let $ Q= \Omega \times (0,1) $, where
$ \Omega = (0,1)^2 \setminus \{(x_1,0)\in \mathbb{R}^2 : 0 \le x_1 < 1\} $
is a ``slit domain'' that is not Lipschitz.
Moreover, we choose the constant diffusion coefficient $ \nu \equiv 1 $, and the manufactured solution
$ u(r,\phi, t) = t\,r^{\alpha}\sin( \alpha \phi),$
where $ (r, \phi) $ are polar coordinates
with respect to
$ (x_1,x_2) $, and $ \alpha = 0.5$.
We know\footnote{\url{https://math.nist.gov/amr-benchmark/index.html}} that
$u$ only belongs to $H^{1+\alpha}(Q)$ due to the singularity of the gradient at $(0,0)$.
Hence, uniform mesh refinement will result in a reduced convergence rate of $ \mathcal{O}(h^{\alpha})$.

\begin{figure}[!htb]
    \includegraphics[width=\linewidth]{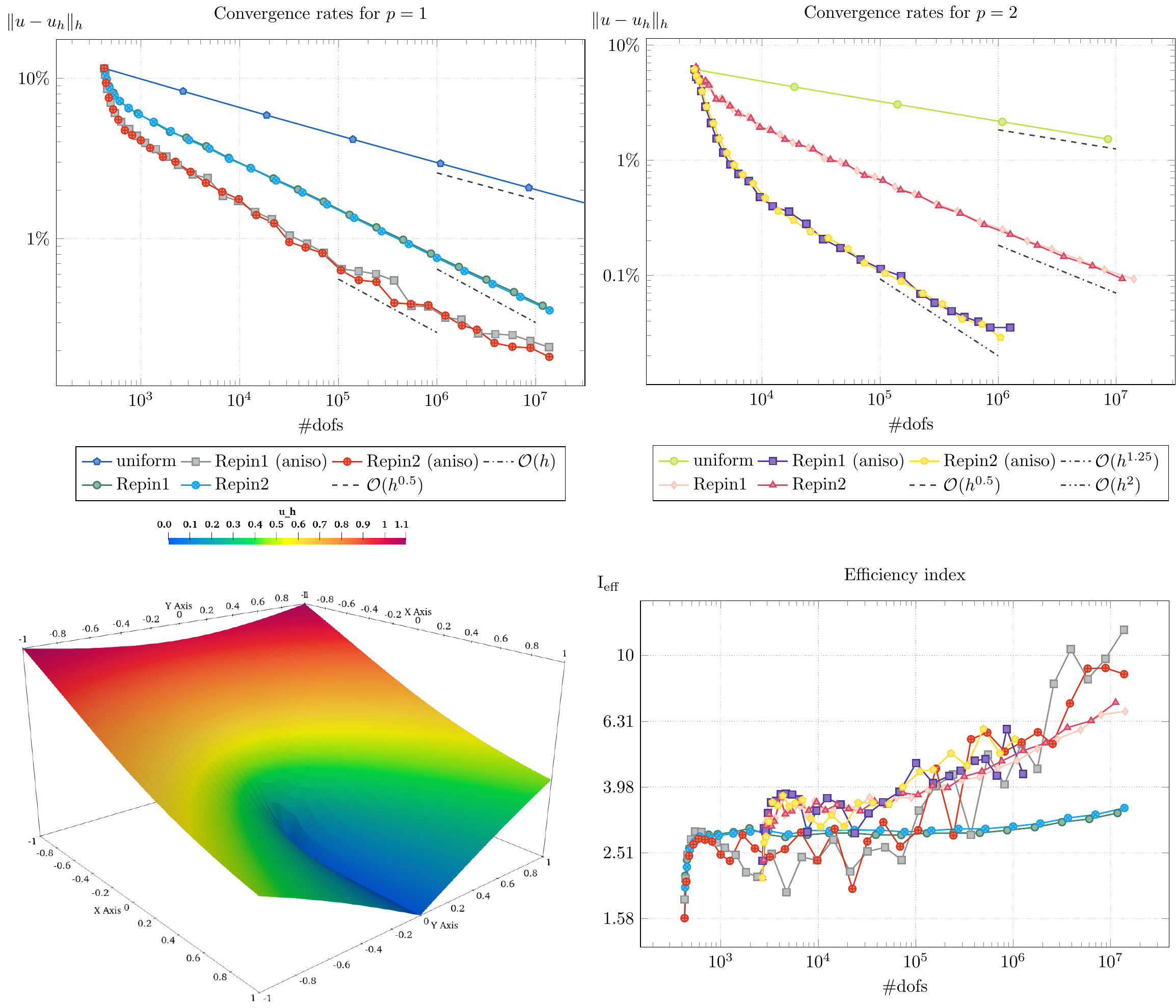}
    \caption{Convergence rates for $ p = 1$ (upper left); convergence rates for $ p = 2 $ (upper right); plot of $ u(\cdot,\cdot, 1) $ (lower left); efficiency indices for $ p= 1,2$, with the respective colors from the upper plots (lower right).}
    \label{LS:fig:1}
\end{figure}

In order to properly realize the adaptive refinement strategies, we need to choose
appropriate $ \mathbf{y} $ and $ \vartheta $. For the first majorant, we reconstruct an improved flux $ \mathbf{y}_h^{(0)} = R_h(\nabla_x u_h) $, where $ R_h $ is a nodal averaging operator. We then
improve this flux by applying a few CG steps to the minimization problem $ \min_{\mathbf{y}} \mathfrak{M}_1 $, obtaining the final flux $ \mathbf{y}_h^{(1)} $ that is then used in the estimator. For the second majorant, we follow the same procedure, but right before postprocssing the flux, we first apply some CG iteration to another minimization problem $ \min_{\vartheta} \mathfrak{M}_2 $.
For linear finite elements, we observe at least optimal convergence rates for both error estimators. Anisotropic refinements, with the anisotropy parameter $ \chi = 0.1 $, manage to obtain a better constant than isotropic refinements; see Fig.~\ref{LS:fig:1} (upper left). For quadratic finite elements, anisotropic adaptive refinements, with $ \chi = 0.15 $, manage to recover the optimal rate of $ \mathcal{O}(h^2) $, while isotropic adaptive refinements result in a reduced rate of $ \mathcal{O}(h^{1.25}) $; see Fig.~\ref{LS:fig:1} (upper right).
The efficiency indices are rather stable for isotropic refinements,
while some oscillations can be observed  for anisotropic refinements; see Fig.~\ref{LS:fig:1} (lower right).
%
%
\bibliographystyle{acm}
\bibliography{ms}
%
\end{document}